\documentclass[11pt]{article}
\usepackage{mathrsfs}
\usepackage{graphicx}
\usepackage{latexsym,amsmath,amssymb,amsfonts}
\usepackage{color,colordvi,amscd,indentfirst,graphics}
\allowdisplaybreaks

\newtheorem{theorem}{Theorem}[section]

\newtheorem{definition}[theorem]{Definition}

\newtheorem{claim}[theorem]{Claim}

\newtheorem{proposition}[theorem]{Proposition}
\newtheorem{problem}[theorem]{Problem}

\def\qed{\hfill \rule{4pt}{7pt}}
\def\pf{\noindent {\it Proof.} }

\begin{document}
	
	\title{Spanning tree with bounded degree sum of two adjacent vertices}
	
	\author{Pham Hoang Ha\footnote{E-mail address: ha.ph@hnue.edu.vn }\\
		Department of Mathematics\\
		Hanoi National University of Education\\
		136 XuanThuy Str., Hanoi, Vietnam}%
	
	\date{}
	\maketitle{}
	\bigskip
	
	\begin{abstract}
		 In this paper, we first give a sharp sufficient condition on degree sum for a connected graph to have a spanning tree whose maximum of degree sum of two adjacent vertices are bounded. These improve some previous results on the spanning tree with bounded maximum degree. Next, we also use our main result to study the generalized rainbow connection of graphs. In particular, we give some degree sum conditions for a connected graphs to have a bounded (1,2)- rainbow connection number.
	\end{abstract}
	
	
	
	\newpage

\section{Introduction}
In this paper, we only consider finite simple graphs. Let $G$ be a graph with vertex set $V(G)$ and edge set $E(G)$. The order of $G$ is denoted by $|G|$. 
For any vertex
$v\in V(G)$, we use $N_G(v)$ and $\deg_G(v)$ (or $N(v)$ and $\deg (v)$ if
there is no ambiguity) to denote the set of neighbors of $v$ and the
degree of $v$ in $G$, respectively. For any $S\subseteq V(G)$, we
denote by $|S|$ the cardinality of $S$. Sometime, we denote by $|G|$ instead of $|V(G)|.$ For a positive integer $k$, we use the notation $[k]$ to represent the set $\{1, 2, \ldots, k\}$. We use $G-S$ to denote the graph obtained from $G$ by deleting the
vertices in $S$ together with their incident edges. We define $G-uv$ to be the
graph obtained from $G$ by deleting the edge $uv\in E(G)$, and
$G+uv$ to be the graph obtained from $G$ by adding a new edge $uv$ joining two non-adjacent vertices $u$ and $v$ of $G$.  For two vertices $u$ and $v$ of $G$, the distance between $u$ and $v$ in $G$ is denoted 
by $d_{G}(u, v)$. 

Let $T$ be a tree. A vertex of degree one is a \emph{leaf} of $T$
and a vertex of degree at least three is a \emph{branch vertex} of
$T$. Setting $L(T)$ and $ B(T)$ the set of leaves and branch vertices of the tree $T$ respectively. For an integer $k\geq 2,$ $k$- tree is a tree with the maximum degree at most $k$. A spanning tree of a graph $G$ is a tree $T$ with $V(T)=V(G).$ In particular, a Hamiltonian path is a spanning $2-$tree. We refer
to~\cite{Di05} for terminology and notation not defined here.

 The maximum independence number of $G$ is denoted by $\alpha(G).$  For each positive integer $p,$ set $\sigma_p(G)$ to be the minimum degree sum of $p$ independent vertices if $\alpha(G)\geq p$ and $\sigma_{p}(G)=+\infty$ otherwise. In particular, $\sigma_{1}(G)=\delta(G),$ the minimum degree of $G.$\\
In 1952, Dirac \cite{D} prove the following theorem.  
\begin{theorem}[{\cite[Dirac]{D}}]\label{Dirac} Let $G$ be a graph of order  $n(n\geq 2)$ and $\delta(G)\geq \dfrac{n-1}{2}.$ Then $G$ has a Hamiltonian path.
	\end{theorem}
Later, Ore~\cite{Ore} obtained a
sufficient condition related to the degree sum for
a connected graph to have a Hamiltonian path. 
\begin{theorem}[{\cite[Ore]{Ore}}]\label{Ore}
	Let $G$ be a connected graph. If $\sigma_2(G) \geq |G|-1$, then $G$ has
	a Hamiltonian path.
\end{theorem}
After that, Win generalized Ore's result by considering a Hamiltonian path as a spanning $2$-tree.
\begin{theorem}[{\cite[Win]{W}}] \label{Win}
 Let $k \geq 2$ and let $G$ be a connected graph of order $n.$ If $\sigma_k(G) \geq n-1$, then $G$ has a spanning $k$-tree.
	\end{theorem}
Regarding to this research direction, many results have been stated (see \cite{CKR1, CKR2, FMMMTY, JW, KOSTY, NR, Yan} for examples).

 On the other hand, a graph $G$ is called an edge-colored graph (or briefly a colored graph) if each of its edge is assigned by one color from $m$ colors, which is represented by the function $c: E(G)\rightarrow [m]$. 

The study of connection problems in graph theory has become an active area of research, attracting considerable attention from mathematicians due to their widespread use in the actual world recently. Now, we present some famous connection problems relating to the motivation of main result in this paper. 

Let $G$ be an edge-colored graph. The initial notion of connection, introduced by Chartrand et al. \cite{Chartrand2008}, is \emph{rainbow connection}. A path within $G$  is called a \emph{rainbow path} if no two edges of it have the same color. We say a connected graph $G$ to be \emph{rainbow connected} if any two distinct vertices are joined by at least one rainbow path. Moreover, the \emph{rainbow connection number} of $G$, denoted by $rc(G)$, is defined as the least number of colors needed to make it rainbow connected. In 2009, Chakraborty et al. in \cite{Chakraborty2009} proved that deciding whether $rc(G) = 2$ is an NP-complete problem. Later, Ananth et al. in \cite{Ananth2011} supplemented the result of Chakraborty et al. and pointed out that the problem remains NP-hard for $rc(G)\geq3$. We refer the readers to references \cite{Li2013, Li2012, Li2017} for surveys in this topic.

Inspired by rainbow connection and proper colorings, the concept of \emph{proper connection} was independently proposed by  Borozan et al. in \cite {Borozan2012} and Andrews et al. in \cite{Andrews2016}. A \emph{proper path} in $G$ is defined as a path $P$ where every pair of consecutive edges are colored differently. A connected graph $G$ in which any two distinct vertices are joined by at least one proper path is called \emph{properly connected}. Furthermore, the \emph{proper connection number} of $G$, denoted by $pc(G)$, is the smallest number of colors required to make $G$ properly connected. Huang et al. in \cite{Huang2020} obtained that calculating $pc(G)$ for a given graph $G$ is recognized as an NP-hard problem. Many results related to this field can be found in \cite{Aardt2017, Li2015, Lumduanhom2016}  for examples.

In 2016, Li et al. in \cite{Li2016} introduced the definition of $(p,l)$-rainbow connection which generalizes both of rainbow connection and proper connection, where $p,l$ are two positive integers. A path $P$ in $G$ is referred to as an \emph{$l$-rainbow path} if any subpath of $P$ containing no more than $l+1$ edges is a rainbow path. A connected graph $G$ is called to be \emph{$(p,l)$-rainbow connected} if for any two distinct vertices, there exist at least $k$ pairwise internally vertex-disjoint $l$-rainbow paths connecting them in $G$. The \emph{$(p,l)$-rainbow connection number} of $G$, denoted by $rc_{p,l}(G)$, is the least number of colors needed to make $G$ $(p,l)$-rainbow connected. The $(1,1)$-rainbow connection number of $G$ is its proper connection number, which derives directly from the definitions above. When $l+1$ is larger than or equal to the diameter of $G$, $(1, l)$-rainbow connection number becomes its rainbow connection number. 

In this research, we focus on the $(1,2)$-rainbow connection number. There are many results of the $(1,2)$-rainbow connection number published in \cite{DDS, Li2016, Li2018, Zhu2021} in recent years. Moreover, in \cite{Li2016}, the authors also proved that if $H$ is a spanning subgraph of $G$ then $ rc_{1,2}(G)\leq rc_{1,2}(H)$ and $(1,2)$-rainbow connection number of a tree $T$ is 1 less than the maximum degree sum of two adjacent vertices in $T.$ This leads to the fact that we can evaluate the upper bound of $rc_{1,2}(G)$ if $G$ has a spanning tree $T$ with bounded degree sum of two adjacent vertices.

Motivate by two above research directions, we would like to study a new result on degree sum condition such that it not only improves previous results on spanning tree with bounded degree but also gives new results of the $(1,2)$- rainbow connection number. In particular, we state a sharp degree sum condition for a graph having a spanning tree with bounded degree sum of two adjacent vertices.  

Before stating the main theorem, we introduce a new notation.
\begin{definition}
	Let $T$ be a tree. We set $\Delta_2(T)=\displaystyle \max_{xy\in E(G)} \{\deg_T(x)+\deg_T(y)\}.$
\end{definition}
It is easy to see that a tree $T$ is a Hamiltonian path if $\Delta_2(T)\leq 4$ unless $T$ is $K_{1,3}.$

The main purpose of this paper is to prove the following result.
\begin{theorem}\label{thm1}
	Let $G$ be a connected graph of order $n$ and $k\geq 2$ be an integer. If $\sigma_k(G)\geq n-1$, then $G$ has a spanning tree $T$ such that
	$$
	\Delta_2(T)\leq k+2.
	$$
\end{theorem}

This paper is organized as follows: In the section 2, we prove Theorem \ref{thm1} and show the sharpness of it. Furthermore, we also point out that it is an improvement of Win's previous result. In the section 3, we focus on studying generalized rainbow connection number of a connected graph for the case $(1,2)$- rainbow connection number by using the main result in the section 2. In the last section, we introduce some comments regarding our results.
\section{Spanning tree with bounded degree sum of two adjacent vertices}
In this section, we focus on studying the conditions of degree sum for a connected graph to have a spanning tree whose degree sum of two arbitrary adjacent vertices are bounded. Namely, we prove Theorem \ref{thm1}. After that, we show that the result of Theorem \ref{thm1} is sharp and it also is an improvement of Theorem \ref{Win}.

To prepare for the proof of main results, we now recall some definitions and notations: 

	Let $T$ be a tree. For two distinct vertices $u,v$ of $T$, we denote by $P_T[u,v]$
	the unique path in $T$ connecting $u$ and $v$ and denote by
	$d_T(u,v)$ the distance between $u$ and $v$ in $T$. We also denote $\{u_v\}=V(P_T[u,v])\cap N_T(u)$. 
	
	We define the
	\emph{orientation} of $P_T[u,v]$ is from $u$ to $v$. If $x \in V(P)$, then $x^+$ and $x^-$ denote the successor and predecessor of $x$ on $P$ if they exist, respectively. For any $X\subseteq V(G),$ set $(N(X)\cap P_T[u,v])^{-} = \{x^{-} \vert x \in V(P_T[u,v])\setminus \{u\}\ \text{and}\ x \in N(X) \}$ and  $(N(X)\cap P_T[u,v])^{+} = \{x^{+} \vert x \in V(P_T[u,v])\setminus \{v\} \ \text{and}\ x \in N(X) \}.$

Let $T$ be a tree and $x\in V(T).$ We set $\Delta_2^T(x)=\displaystyle \max_{y\in N_T(x)} \{\deg_T(x)+\deg_T(y)\}.$ Then  $\Delta_2(T)=\displaystyle \max_{x\in V(T)} \{\Delta_2^T(x)\}.$

\subsection{Proof of Theorem \ref{thm1}}
	
We prove the theorem
by contradiction. Suppose to the contrary that $G$ contains no
spanning tree $S$ with $\Delta_{2}(S)\leq k+2.$ 

Now, we may choose a maximal tree $T$ of $G$ so that $\Delta_2(T)\leq k+2$ satisfying
\begin{itemize}
\item[$($C1$)$] $|L(T)|$ is as small as possible.
\end{itemize}

Since $T$ is not a spanning tree of $G$, there exist two vertices $y\in V(T)$ and $w\not\in V(T)$ such that $yw \in E(G).$ \\
Now, if $\Delta_2^T(y) < k+2$ then the tree $T'=T+yw$ satisfies that $\Delta_2(T')\leq k+2$ and $|V(T')|> |V(T)|.$ This contradicts to the maximality of $T.$ Then, we obtain $\Delta_2^T(y) = k+2.$

We consider following two cases:

{\bf Case 1.}\ There exists a leaf $u$ of $T$ such that $\Delta_2^T(u) = k+2.$ Then, $T$ must be a star. Hence $T=K_{1,k+1}$ with branch vertex $b$ and $k+1$ leaves $\{u_1,...,u_{k+1}\}.$

{\it Subcase 1.1.}\ $y=b.$

Assume that there exists some edge $u_iu_j \in E(G),$ where $i,j \in [k+1]$. Consider the tree $T'=T-\{u_iy\} + \{u_iu_j\}$ then $|V(T')|= |V(T)|, \Delta_2(T')=\Delta_2(T)$ and $|L(T')|<|L(T)|,$ a contradiction with the condition (C1). This implies that $u_iu_j \not\in E(G)$ for every $i, j \in [k+1].$ Namely, $\{u_1,...,u_{k+1}\}$ is an independent set.

Moreover, if we assume that there is a vertex $v\in V(G)\setminus V(T)$ such that $|N_G(v)\cap L(T)|\geq 2,$ then there exist $1\leq i< j\leq k+1$ satisfying $u_iv, u_jv\in E(G).$ Hence, the tree 
$$T'= T+\{u_iv,u_jv\}-\{u_jy\}$$
 has $\Delta_2(T') = k+2$ and $|V(T')|> |V(T)|.$ This contradicts the maximality of $T.$ \\
 Thus, $|N_G(v)\cap L(T)|\leq 1$ for every $v\in V(G)\setminus V(T).$ In particular, $|N_G(w)\cap L(T)|\leq 1$. Hence, we may assume that $u_iw\not\in E(G)$ for every $2\leq i\leq k+1.$ So we conclude that 
 $$\displaystyle\sum_{i=2}^{k+1}\bigg|N_G(u_i)\cap \bigg(V(G)\setminus V(T)\bigg)\bigg| \leq |V(G)|-|V(T)|-1.$$

Therefore, we obtain

\begin{align*}
	\sigma_{k}(G)&\leq \sum_{i=2}^{k+1}|N_G(u_i)|\\
	&= \sum_{i=2}^{k+1}|N_G(u_i)\cap V(T)|+\sum_{i=2}^{k+1}|N_G(u_i)\cap \bigg(V(G)\setminus V(T)\bigg)|\\
	&\leq k + (|V(G)|-|V(T)|-1) = n- 3.
	\end{align*}
This violates to the assumption of Theorem \ref{thm1}.

{\it Subcase 1.2.}\ $y=u_t$ for some $t\in [k+1].$ For the sake of convenience, we may assume by symmetry that $t=k+1.$

Suppose that there exist indices $1\leq i< j \leq k+1$ such that  $u_iu_j \in E(G),$ we consider the tree 
$$T'= T+\{u_iu_j, yw\}- \{u_jb\}.$$
Then $T'$ has $\Delta_2(T')\leq k+2$ and $|V(T')|> |V(T)|,$ a contradiction. Then we obtain $u_iu_j \not\in E(G)$ for every $1\leq i< j \leq k+1.$ Namely, the set $\{u_1,...,u_{k+1}\}$ is independent.

Moreover, if we assume that there is a vertex $v\in V(G)\setminus V(T)$ such that $|N_G(v)\cap L(T)|\geq 2$ then there exist $1\leq i< j\leq k+1$ satisfying $u_iv, u_jv\in E(G).$ Hence, the tree 
$$T'= T+\{u_iv,u_jv\}-\{u_ib\}$$
 has $\Delta_2(T') \leq k+2$ and $|V(T')|> |V(T)|.$ This contradicts to the maximality of $T.$ So we conclude that $$ \bigg(N_G(u_i)\cap (V(G)\setminus V(T))\bigg)\cap \bigg(N_G(u_j)\cap (V(G)\setminus V(T))\bigg)=\emptyset \text{ for every } 1\leq i< j\leq k+1.$$
So we have $u_iw\not\in E(G)$ for every $1\leq i\leq k.$
Therefore, we obtain

\begin{align*}
	\sigma_{k}(G)&\leq \sum_{i=1}^{k}|N_G(u_i)|\\
	&= \sum_{i=1}^{k}|N_G(u_i)\cap V(T)|+\sum_{i=1}^{k}|N_G(u_i)\cap \bigg(V(G)\setminus V(T)\bigg)|\\
	&\leq k + (|V(G)|-|V(T)|-1) \leq n- 3.
\end{align*}
This gives a contradition with the assumption of Theorem \ref{thm1}.

{\bf Case 2.}\ There is not any leaf $u$ of $T$ such that $\Delta_2^T(u) = k+2.$ Then we have $\Delta_2^T(u) \leq k+1$ for every $u\in L(T).$
\begin{claim}\label{claim1} 
	$L(T)$ is independent and $N_{G}(L(T))\cap (V(G)\setminus V(T))=\emptyset$.
\end{claim}
\pf
Assume that there exist two leaves $u$ and $v$ that are adjacent in $G$. If $T$ is a path then we have $T'=T+yw$ is also a path with $|T'|> |T|.$ This gives a contradiction to the maximality of $T.$ Otherwise, then $u$ has some nearest branch vertex $b$.
Let $T'=T-\{bb_u\}+\{uv\}.$ Then $T'$ satisfies $\Delta_2(T')\leq k+2.$ Moreover, we also obtain that $|L(T')|<|L(T)|,$ the reason is that either $T'$ has only one new leaf $b_u$ and $u, v$ are not leaves of $T'$ or $u$ is still a leaf of $T'$ but $T'$ has no new leaf and $v$ is not a leaf of $T'$. This contradicts to the condition (C1). Hence $L(T)$ is independent. 

Beside that, if there exist two vertices $u\in L(T)$ and $v\in V(G)\setminus V(T)$ such that $uv\in E(G).$ Then the tree $T'=T+\{uv\}$ has $\Delta_2(T')\leq k+2$ by the assumption of \textbf{Case 2}. This contradicts to the maximality of $T.$ Therefore, we conclude $N_{G}(L(T))\cap (V(G)\setminus V(T))=\emptyset.$ 

So the claim holds.
\qed

Let $\{x_i\}_{i=1}^m$ be the subset of $N_T(y)$ such that $\deg_T(y)+\deg_T(x_i)=k+2$ (we note that $m$ depend on $y$). Thus $\Delta_2^T(x_i)=k+2$ . Then there exists a positive integer $q$ such that $\deg_T(x_i)=q+1$ for every $ i\in [m].$ For each vertex $x_i, i \in [m]$,  set $B_{ij} (1\leq j \leq k)$ to be components of $T[V(T)\setminus \{x_i,y\}].$ Set $v_{ij}\in V(B_{ij})\cap (N_T(y)\cup N_T(x_i)).$ For each $i\in [m],$ take $u_{ij}\in V(B_{ij})\cap L(T)$ such that $\displaystyle\sum_{j=1}^kd_T(y,u_{ij})$ is as large as possible. Let 

$$l_{(T,y)}=\max_{x_i, i \in [m]}\sum_{j=1}^kd_T(y,u_{ij}),$$

for every $y\in V(T)$ such that $N_G(y)\cap \bigg(V(G)\setminus V(T)\bigg) \not= \emptyset.$

Then, we set 
$$l_{T}=\max_{y\in V(T), N_G(y)\cap \bigg(V(G)\setminus V(T)\bigg) \not= \emptyset }l_{(T,y)}.$$
Denote by $B_{ij}, 1\leq j\leq q$ the components which are incident to $x_i$ in $T.$ 

Now, we choose $(T, y)$ such that
\begin{itemize}
		\item[$($C2$)$] $l_{T}=l_{(T,y)}$ is as large as possible, subject to (C1).
\end{itemize}

\begin{claim}\label{claim2} There exists some vertex $x_i, i\in [m]$ such that, for every $t, j\in [q]$ and $ t\not= j,$ we have
	$$N_G(u_{ij})\cap \bigg(V(B_{it})\setminus \{u\in V(B_{it})\, |\, \Delta^T_2(u)=k+2\}\bigg)=\emptyset.$$
\end{claim}
\pf
Suppose to the contrary that for every $i\in [m],$ there exist a vertex $u_{ij}$ and a vertex $a_{it}\in V(B_{it})$ such that $t\not= j, \Delta^T_2(a_{it})<k+2$ and $u_{ij}a_{it}\in E(G),$ (we note that $a_{it}\not\in \{x_i\}_{i=1}^m$). Now we consider the tree 
 $$T'=T+ \{\displaystyle\cup_{i\in [m]}u_{ij}a_{it}\} +\{yw\} - \{\displaystyle\cup_{i\in [m]}x_iv_{ij}\}.$$

Then $|\Delta_2(T')|\leq k+2$ and $|V(T')|=|V(T)|+1,$ the reason is that $\deg_T(y)+\deg_T(c) < k+2$ for every $c\in N_T(y)\setminus \{x_i\}_{i=1}^m.$ This gives a contradiction. So the assertion of Claim \ref{claim2} holds.
\qed

Without loss of generality, let $x_1$ be such a vertex. We now set $P_{1j}=P_T[v_{1j},u_{1j}]$ for each $j \in [k].$ 
\begin{claim}\label{claim3} For every $u_{1j}, j\in [k]$ and $z\in V(T)\setminus \bigg( \displaystyle\cup_{t=1}^{k}V(P_{1t}\cup \{x_1,y\})\bigg),$ we have $u_{1j}z\not\in E(G).$
	\end{claim}
\pf
	Suppose for a contradiction that there exist some vertices $u_{1j}, j\in [k]$ and $z\in V(T)\setminus \bigg( \displaystyle\cup_{t=1}^{k}V(P_{1t})\cup \{x_1,y\})\bigg)$ such that $u_{1j}z \in E(G).$ Then $T+\{u_{1j}z\}$ contains a unique cycle $C$. Let $e$ be an edge of $C$ incident with $z$ and $e\not= u_{1j}z.$ Consider the tree $T'=T+\{u_{1j}z\}-\{e\}.$ Then $\Delta_2(T')\leq k+2,$ $|L(T')|\leq |L(T)|$ and $l_{T'}> l_{T}.$
	 This contradicts the conditions either (C1) or (C2). Therefore, this completes the proof of Claim \ref{claim3}.
\qed

By Claim \ref{claim2}, we have the following claim:
\begin{claim}\label{claim4} For every vertices $u_{1j}, j\in [k]$ and $z\in V(P_{1t}), t\in [k]\setminus \{j\},$ if $u_{1j}z\in E(G),$ then $\Delta_2^T(z)=k+2.$
	\end{claim}
\begin{claim}\label{claim5} We have $\sum_{j=1}^k|N_G(u_{1j})\cap \{y,x_1\}|\leq k+1.$ In special case, if there exist two indices $s \in [q]$ and $t\in [k]\setminus [q]$ such that $u_{1t}v_{1s}\in E(G),$ then $\sum_{j=1}^k|N_G(u_{1j})\cap \{y,x_1\}|\leq k.$ 
\end{claim}

\pf

Suppose that there exist $u_{1j}y, u_{1t}y\in E(G)$ for some distinct indices $j, t\in [q].$ Then $\deg_T(x_1)\geq 3.$ Consider the tree $T':=T+\{u_{1j}y\}-\{x_1y\}$ has $\Delta_2(T')\leq k+2$ and $|L(T')|< |L(T)|,$ a contradiction. So we obtain $\displaystyle\sum_{j=1}^q|N_G(u_{1j})\cap \{y\}|\leq 1.$ Hence, 
$$\displaystyle\sum_{j=1}^k|N_G(u_{1j})\cap \{y\}|\leq k-q+ 1.$$

On the other hand, suppose that there exists $u_{1j}x_1\in E(G)$ for some $j\in [k]\setminus [q].$ Then the tree $T'=T+\{u_{1j}x_1, yw\}-\{x_1y\}$ has $\Delta_2(T')\leq k+2$ and $|V(T')|> |V(T)|.$ This contradicts to the maximality of $T.$ 
So $u_{1j}x_1\not\in E(G)$ for every $j\in [k]\setminus [q].$ Hence we have $$\sum_{j=1}^k|N_G(u_{1j})\cap \{x_1\}|= \sum_{j=1}^q|N_G(u_{1j})\cap \{x_1\}| \leq q.$$ 
Moreover, if there exist two indices $s \in [q]$ and $t\in [k]\setminus [q]$ such that $u_{1t}v_{1s}\in E(G),$ then we obtain that $x_1u_{1s}\not\in E(G).$ Indeed, suppose that the assertion is false. Consider the tree $T'=T+\{x_1u_{1s}, u_{1t}v_{1s}, yw \}-\{x_1v_{1s}, x_1y\}.$ Then $T'$ has $\Delta_2(T')\leq k+2$ and $|V(T')|> |V(T)|.$ This contradicts to the maximality of $T.$ Hence $x_1u_{1s}\not\in E(G).$ Therefore, we have $$\sum_{j=1}^k|N_G(u_{1j})\cap \{x_1\}|= \sum_{j=1}^q|N_G(u_{1j})\cap \{x_1\}| \leq q-1.$$ 

Combining two above conclusions, we get $\displaystyle\sum_{j=1}^k|N_G(u_{1j})\cap \{y,x_1\}|\leq k+1$\\ and  $\displaystyle\sum_{j=1}^k|N_G(u_{1j})\cap \{y,x_1\}|\leq k$ if there exist two indices $s \in [q]$ and $t\in [k]\setminus [q]$ such that $u_{1t}v_{1s}\in E(G).$ 

 Therefore, Claim \ref{claim5} holds.
\qed
\begin{claim}\label{claim6} For each $j\in [k],$ if there exist vertices $u_{1i}, i \in [k]\setminus \{j\}$ and $z\in V(P_{1j})$ such that $u_{1i}z\in E(G),$ then $u_{1t}z^-\not\in E(G)$ and $u_{1t}z^+ \not\in E(G)$ for all $t\not=i, t\not=j.$ Moreover, $\deg_T(z^-)= 2.$
\end{claim}
\pf
By symmetry, we only need to prove that $u_{1t}z^-\not\in E(G)$ for some $t\in [k]\setminus \{i, j\}$ Indeed, we consider two cases:\\
If $j\in [q].$ Suppose that assertion is false. Then $u_{1j}z^-\in E(G).$ If $\deg_T(x_1) \geq 3,$ consider the tree $T'=T+\{u_{1j}z^-, u_{1i}z\}-\{zz^-, x_1v_{1j}\}.$ Then $T'$ has  $\Delta_2(T')\leq k+2$ and $|L(T')|< |L(T)|,$  which contradicts (C1). Otherwise, $\deg_T(x_1) =2,$ then we consider the tree $T'=T+\{u_{1j}z^-, u_{1i}z, yw\}-\{zz^-, yx_1\}.$ Then $T'$ has  $\Delta_2(T')\leq k+2$ and $|V(T')|> |V(T)|,$  a contradiction.\\
If $j\in [k]\setminus [q].$ Suppose that contrary that $u_{1j}z^-\in E(G).$ Consider the tree $T'=T+\{u_{1j}z^-, u_{1i}z,yw\}-\{zz^-, yv_{1j}\}.$ Then $T'$ has  $\Delta_2(T')\leq k+2$ and $|V(T')|> |V(T)|,$ a contradiction.\\
So the first part of claim holds.

Now, if $\deg_T(z^-) \geq 3$ then the tree $T'=T+\{u_{1i}z\}-\{zz^-\}.$ Then $T'$ has  $\Delta_2(T')\leq k+2$ and $|L(T')|< |L(T)|,$  contradicting with (C1). 

Claim \ref{claim6} is proved. 
\qed

\begin{claim}\label{claim7} For each $j\in [k],$ if $u_{1j}z\in E(G)$ for some $z\in V(P_{1j})$ then $u_{1i}z^+ \not\in E(G)$ for every vertex $u_{1i}, i\not=j$ and $\deg_T(z^+)=2.$ 
\end{claim}
\pf
Suppose to the contrary that $u_{1i}z^+\in E(G).$ \\
If $j\in [q],$ then the tree 
$$T'= \left\{\begin{array}{ll} T+\{u_{1i}z^+, u_{1j}z\}-\{zz^+, x_1v_{1j}\},
	& \;\mbox{ if } i\in [q],\\
	T+\{u_{1i}z^+, u_{1j}z, yw\}-\{zz^+, yx_{1}\},
	& \;\mbox{ if } i\in[k]\setminus [q],
\end{array}\right.$$
 satisfies  $\Delta_2(T')\leq k+2.$ Moreover, we have $|L(T')|< |L(T)|$ when $i\in [q]$ and $|V(T')|> |V(T)|$ when $i\in[k]\setminus [q].$ This contradicts either (C1) or the maximality of $T$. \\
 If $j\in [k]\setminus [q],$ then the tree 
 $$T'= \left\{\begin{array}{ll} T+\{u_{1i}z^+, u_{1j}z, yw\}-\{zz^+, yv_{1j}\},
 	& \;\mbox{ if } i \in [k]\setminus [q],\\
 	T+\{u_{1i}z^+, u_{1j}z, yw\}-\{zz^+, yx_{1}\},
 	& \;\mbox{ if } i\in [q],
 \end{array}\right.$$
 satisfies  $\Delta_2(T')\leq k+2.$ Moreover, we have $|V(T')|> |V(T)|$. This contradicts to the maximality of $T$.

Now, if $\deg_T(z^+) \geq 3,$ then the tree $T'=T+\{u_{1j}z\}-\{zz^+\}.$ Hence $T'$ has  $\Delta_2(T')\leq k+2$ and $|L(T')|< |L(T)|,$  a contradiction with (C1). 

Therefore, we complete the proof of Claim \ref{claim7}. 

\qed

Set $U=\{u_{1i}\}_{i=1}^k.$ \\
Fix $j\in [k].$ Set $V(P_{1j})=\{v_1, v_2, ..., v_h\},$ where $v_1=v_{1j},  v_h= u_{1j}$ and $v_i = v^-_{i+1},$ for every $1\leq i\leq h-1.$ Let $A=\{a\in V(P_{1j})| |N_G(a)\cap U| \geq 2\}=\{v_{j_1}, ..., v_{j_r}\}.$\\
By Claim \ref{claim1} we have $j_r < h.$

\begin{claim}\label{claim8} We have
\begin{equation*}
		\displaystyle\sum_{u\in U}|N_G(u)\cap V(B_{1j})| \leq |V(B_{1j})|-1, \,\, \text{ for all } j\in [k]\setminus [q].
	\end{equation*}
\end{claim}
\pf
Assume that $u_{1t}v_{1j}\in E(G)$ for some $t\not=j.$ Then the tree $ T'= T+\{u_{1t}v_{1j}, yw\}-\{yv_{1j}\}$
satisfies  $\Delta_2(T')\leq k+2$ and $|V(T')|> |V(T)|,$  which contradicts to the maximality of $T$. Hence, we conclude that  $u_{1t}v_{1j}\not\in E(G)$ for every $t\not=j.$ This implies that $j_1 > 1.$

Since Claim \ref{claim6} holds, we conclude that $j_{t+1}-j_t > 1$ for every $t\leq r-1,$ in particular, if $v_i \in A$ then $v_{i+1}\not\in A$ for all $i\in [h-1].$

Take a vertex $v_{j_i}\in A.$ By Claim \ref{claim4} we have $\Delta_2^T(v_{j_i})=k+2.$ We now consider following three facts:\\
{\it Fact 1.}\, If $\deg_G(v_{j_i-1})+\deg_{G}(v_{j_i})=k+2,$ set the pair of vertices $(v_{j_i-1}, v_{j_i})$ in $A_1.$\\
Then $v_{j_i-1}\not\in \cup_{i\in [k]}N_G(u_{1i})$ by Claims \ref{claim6} and \ref{claim7}. Hence we obtain
\begin{align*}
	\sum_{t=1}^k|N_G(u_{1t})\cap \{v_{j_i-1}, v_{j_i}\}|&=\sum_{t=1}^k|N_G(u_{1t})\cap \{v_{j_i}\}|\\
	&\leq k = \bigg|\bigg(( N_{T}(v_{j_i-1})\cup N_{T}(v_{j_i}))\setminus V(P_{1j})\bigg)\cup \{v_{j_i-1}, v_{j_i}\}\bigg|.
	\end{align*}
{\it Fact 2.}\, If $\deg_{G}(v_{j_i-1})+\deg_{G}(v_{j_i})\not= k+2$ and $\deg_{G}(v_{j_i+1})+\deg_{G}(v_{j_i})= k+2,$ set the three vertices $(v_{j_i-1}, v_{j_i}, v_{j_i+1})$ in $A_2.$\\
By Claims \ref{claim6} and \ref{claim7}, we obtain $\displaystyle\sum_{t=1}^k|N_G(u_{1t})\cap \{v_{j_i-1}, v_{j_i}\}|\leq k$ and $v_{j_i+1}\not\in \cup_{i\in [k]\setminus \{j\}}N_G(u_{1i})$. Thus we have
\begin{align*}\displaystyle\sum_{t=1}^k|N_G(u_{1t})\cap \{v_{j_i-1}, v_{j_i}, v_{j_i+1}\}|= \displaystyle\sum_{t=1}^k|N_G(u_{1t})\cap \{v_{j_i-1}, v_{j_i}\}|+\displaystyle\sum_{t=1}^k|N_G(u_{1t})\cap \{v_{j_i+1}\}|\\
	\leq k+1 = \bigg|\bigg( (N_{T}(v_{j_i+1})\cup N_{T}(v_{j_i}))\setminus V(P_{1j})\bigg)\cup \{v_{j_i-1}, v_{j_i}, v_{j_i+1}\}\bigg|.\end{align*}
{\it Fact 3.}\, If $\deg_{G}(v_{j_i-1})+\deg_{G}(v_{j_i})\not= k+2$ and $\deg_{G}(v_{j_i+1})+\deg_{G}(v_{j_i})\not= k+2,$ set the vertex $ v_{j_i}$ in $A_3.$\\
Then there exists a vertex $u\in N_T(v_{j_i})\cap \bigg(V(B_{1j})\setminus V(P_{1j})\bigg)$ such that $\deg_G(u)+\deg_G(v_{j_i})= k+2.$ Thus we obtain
$$\sum_{t=1}^k|N_G(u_{1t})\cap \{ v_{j_i}\}|\leq k \leq \bigg|\bigg(( N_{T}(u)\cup N_{T}(v_{j_i}))\setminus V(P_{1j})\bigg)\cup \{v_{j_i}\}\bigg|.$$

On the other hand, if $v_{j_i}, v_{j_i+2}\in A$ and $\Delta_2^T(v_{j_i+1})=k+2$ then $\deg_T(v_{j_i+1})=2$ by Claim \ref{claim7}. Combining with $\Delta_2(T)\leq k+2,$ we have $\deg_T(v_{j_i-1})=2.$ Hence, we consider $(v_{j_i-1}, v_{j_i}),(v_{j_i+1}, v_{j_i+2}) $ as two pair of vertices in {\it Fact 1}.
Therefore, we conclude that $V(P_{1j})$ can be divided into discrete subsets consisting of vertices belonging to $A_3$, pairs of vertices belonging to $A_1$, sets of three vertices belonging to $A_2$ and $V(P_{1j})\setminus (A_1\cup A_2\cup A_3).$

By definition, we have $|N_G(v)\cap U|\leq 1$ for all $v\in V(P_{1j})\setminus (A_1\cup A_2\cup A_3)$ and $N_G(u_{1j})\cap U = \emptyset.$

Combining with {\it Facts 1-3} and Claim \ref{claim3}, we now obtain that
\begin{align*}
	&\displaystyle\sum_{u\in U}|N_G(u)\cap V(B_{1j})|= \displaystyle\sum_{u\in U}|N_G(u)\cap V(P_{1j})| \\
	&\leq \displaystyle\sum_{(v_{j_i-1}, v_{j_i})\in A_1}|N_G(u)\cap \{v_{j_i-1}, v_{j_i}\}|+\displaystyle\sum_{(v_{j_i-1}, v_{j_i}, v_{j_i+1})\in A_2}|N_G(u)\cap \{v_{j_i-1}, v_{j_i}, v_{j_i+1}\}|\\ 
	&+\displaystyle\sum_{v_{j_i}\in A_3}|N_G(u)\cap \{v_{j_i}\}|+|V(P_{1j})\setminus(A_1\cup A_2\cup A_3)| \\
		&\leq |V(B_{1j})\setminus \{u_{1j}\}|=|V(B_{ij})|-1. 
\end{align*}
Therefore, we conclude that
\begin{equation*}\label{eq1}
	\displaystyle\sum_{u\in U}|N_G(u)\cap V(B_{1j})| \leq |V(B_{1j})|-1, \,\, \text{ for all } j\in [k]\setminus [q].
\end{equation*} 
Claim \ref{claim8} holds.
\qed
\begin{claim}\label{claim9} If $\deg_T(x_1)\geq 3$ or $v_{1j}\not\in A$ then	\begin{equation*}
		\displaystyle\sum_{u\in U}|N_G(u)\cap V(B_{1j})| \leq |V(B_{1j})|-1, \,\, \text{ for all } j\in [q].
	\end{equation*}
\end{claim}
\pf
 Suppose that $u_{1t}v_{1j}\in E(G)$ for some $t\not=j$ when $\deg_T(x_1)\geq 3.$ Consider the tree $ T'= T+\{u_{1t}v_{1j}\}-\{x_1v_{1j}\}.$ Then $T$
has  $\Delta_2(T')\leq k+2$ and $|L(T')|< |L(T)|,$  which contradicts to (C1). Hence, we obtain that  $u_{1t}v_{1j}\not\in E(G)$ for every $t\not=j.$ This implies that $j_1 > 1.$ On the other hand, if $v_{1j}\not\in A$ then $j_1 > 1.$ So we always conclude $j_1 > 1.$  Now, repeating the arguments as in the proof of Claim \ref{claim8} we complete the proof of Claim \ref{claim9}. 
\qed
\begin{claim}\label{claim10} If $\deg_T(x_1)= 2$ and $v_{11}\in A$ then $q=1$ and 	\begin{equation*}
		\displaystyle\sum_{u\in U}|N_G(u)\cap V(B_{11})| \leq |V(B_{11})|.
	\end{equation*}
\end{claim}
\pf
It is easy to see that $q=1$ from $\deg_T(x_1)= 2.$\\
Suppose that there exists some $t\in [k]\setminus\{1\}$ such that $u_{1t}x_1\in E(G).$ Consider the tree $ T'= T+\{u_{1t}x_1, yw\}-\{yx_1\}.$ Then $T'$ satisfies  $\Delta_2(T')\leq k+2$ and $|V(T')|> |V(T)|.$ This contradicts to the maximality of $T.$ So $u_{1t}x \notin E(G)$ for every $t\not= 1.$\\
Set  $\overline{P}_{11}$ = $P_{11}+\{x_1v_{11}\}.$ Repeating the arguments as in the proof of Claim \ref{claim8} when we replace $P_{11}$ by $\overline{P}_{11},$ we obtain that 
\begin{equation*}
	\displaystyle\sum_{u\in U}|N_G(u)\cap V(B_{11}+\{x_1v_{11}\})| \leq |V(B_{11}+\{x_1v_{11}\})|-1 = |V(B_{11})|.
\end{equation*}
Therefore, the proof of Claim \ref{claim10} is finished.
\qed

Next we combine all above claims to complete the proof of Theorem \ref{thm1}.\\
If $\deg_T(x_1)\geq 3$ or $v_{1j}\not\in A,$ then by Claims \ref{claim1}, \ref{claim5}, \ref{claim8}, \ref{claim9} we get
	\begin{align*}
		\sigma_{k}(G) & \leq \displaystyle\sum_{u\in U}|N_G(u)|= \displaystyle\sum_{t=1}^k\bigg(\sum_{j=1}^k|N_G(u_{1j})\cap V(B_{1t})|\bigg) +\sum_{j=1}^k|N_G(u_{1j})\cap \{y,x_1\}|\\
			&\leq \sum_{t=1}^k(|V(B_{1t})|-1 ) + k+1=\sum_{t=1}^k|V(B_{1t})|+1\\
			& = |T|-1\leq |G|-2 =n-2.
		\end{align*}
	If $\deg_T(x_1)=2$ and $v_{11}\in A,$ then there exists $t\not=1$ such that $u_{1t}v_{11}\in E(G).$  By Claims \ref{claim1}, \ref{claim5}, \ref{claim8}, \ref{claim10} we get
	\begin{align*}
		\sigma_{k}(G) & \leq \displaystyle\sum_{u\in U}|N_G(u)|= \displaystyle\sum_{t=1}^k\bigg(\sum_{j=1}^k|N_G(u_{1j})\cap V(B_{1t})|\bigg) +\sum_{j=1}^k|N_G(u_{1j})\cap \{y,x_1\}|\\
		& \leq \displaystyle\sum_{j=1}^k|N_G(u_{1j})\cap V(B_{11})|+ \displaystyle\sum_{t=2}^k\bigg(\sum_{j=1}^k|N_G(u_{1j})\cap V(B_{1t})|\bigg) +\sum_{j=1}^k|N_G(u_{1j})\cap \{y,x_1\}|\\
		&\leq |V(B_{11})|+\sum_{t=2}^k(|V(B_{1t})|-1 ) + k=\sum_{t=1}^k|V(B_{1t})|+1\\
		& = |V(T)|-1\leq |G|-2 =n-2.
	\end{align*}
This violates the hypothesis of Theorem \ref{thm1}.
 So we complete the proof of Theorem \ref{thm1}.
\qed
\subsection{Sharpness of Theorem \ref{thm1}}
We now show that the condition $\sigma_k(G)\geq |G|- 1$ of Theorem \ref{thm1} is sharp. In particular, we prove the following proposition.
\begin{proposition}\label{p1} Let $k\geq 2$ be an integer. There exist many graphs $G$ such that $\sigma_k(G)= |G|- 2$ and $G$ has no spanning tree $T$ satisfying $\Delta_2(T)\leq k+2.$
	\end{proposition}
\pf Let $m$ be an integer such that $m\geq k\geq 2.$ Set $G=K_{m,m(k-1)+2}$ to be a complete bipartite graph with bipartition $V(K_{m,m(k-1)+2})=A\cup B, |A|=m, |B|=m(k-1)+2.$ Then $n=|G|=km+2$ and $\sigma_k(G)=km= |G|- 2.$ Consider each spanning tree $T$ of $G.$ Since $T$ is a tree, we have $|E(T)|=|T|-1 = km+1$ and 
 $$2|E(T)|=\displaystyle\sum_{x\in V(T)}\deg_T(x)=\displaystyle\sum_{y\in A}\deg_T(y)+\displaystyle\sum_{z\in B}\deg_T(z).$$
 On the other hand, since $G$ is bipartite. we also have 
  $$\displaystyle\sum_{y\in A}\deg_T(y)=\displaystyle\sum_{z\in B}\deg_T(z).$$
  Hence, we obtain 
  $$\displaystyle\sum_{y\in A}\deg_T(y)= |E(T)|=km+1.$$
  Therefore, there exists some vertex $y\in A\subset V(T)$ such that $\deg_T(y)\geq k+1.$\\
  
  Now, if $\deg_T(x)=1$ for every $x\in N_T(y),$ then $T=K_{1,d}$, where $d=\deg_T(y).$ Combining with $G$ is bipartite and $\sigma_k(G)= |G|- 2$, we obtain $T=K_{1,k+1}$. This contradicts with $|V(T)|=km+2 > k+2.$ So there exists some vertex $x\in N_T(y)$ such that $\deg_T(x)\geq 2.$ This implies that $\deg_T(x)+\deg_T(y)\geq k+1+2=k+3.$ Hence, we conclude that  $\Delta_2(T)\geq k+3.$ Therefore, $G$ has no spanning tree $T$ such that $\Delta_2(T)\geq k+2.$ The proposition holds. 
  \qed
 \subsection{Theorem \ref{thm1} is an improvemnet of Theorem \ref{Win}} 
   We first recall Theorem \ref{Win} of Win and prove it by using the result of Theorem \ref{thm1}.
\begin{theorem}[{\cite[Win]{W}}] \label{W}
	Let $k \geq 2$ and let $G$ be a connected graph of order $n.$ If $\sigma_k(G) \geq n-1$, then $G$ has a spanning $k-$tree.
\end{theorem}

\pf   Suppose the assertion of Theorem \ref{W} is false. Then, by Theorem \ref{thm1}, there exists some spanning $T$ tree of $G$ such that $\Delta_2(T)\leq k+2$ and $\Delta (T)\geq k+1.$ Since $\Delta_2(T)\leq k+2$ and $\Delta (T)\geq k+1,$ there exists some vertex $a\in V(T)$ such that $\deg_T(a) = k+1.$ Since $\Delta_2(T)\leq k+2$, we obtain $\deg_T(b)=1$ for every $b\in N_T(a).$  Then $T=K_{1,k+1}$ with $L(T)=\{b_1,...,b_{k+1}\}$. Now, if $b_ib_j\in E(G)$ for some $1\leq i< j\leq k+1$ then the tree $T':= T +\{b_ib_j\}-\{ab_i\}$ satisfying  $\Delta (T')\leq k+1,$ a contradiction. So  $b_ib_j\not\in E(G)$ for every $1\leq i< j\leq k+1.$ This implies that  $\sigma_k(G)=k=n-2,$ a contradiction. Therefore, we complete the proof of Theorem \ref{W}.
  \qed
   
 \section{(1,2)-rainbow connection}

 Let $G$ be a connected graph with $n$ vertices and $m$ edges. By the definitions, the following apparent bounds are obtained:
 $$1 \leq rc_{1,1}(G)=pc(G) \leq rc_{1,2}(G)\leq\cdots\leq rc_{1,l}(G)\leq \cdots \leq rc(G) \leq m.$$
 
 Regarding the $(1,2)$-rainbow connection number, many results have been introduced. The following result furnishes the $(1,l)$-rainbow connection number for traceable graphs.
 \begin{proposition}[\cite{Li2016}] \label{pptraceable}
 	If $G$ is traceable, it means that $G$ has a Hamiltonion path, then $rc_{1,l}(G) \leq l+1$ for every positive integer $l.$ Particularly, $rc_{1,2}(G) \leq 3.$
  \end{proposition}
 For the case of tree, we have the following.
 \begin{theorem}[\cite{Li2016}] 
 	\label{theorem_tree}
 	If $T$ is a nontrivial tree, then $rc_{1,2}(T) = \Delta_2(T)-1$. 
 \end{theorem}
The relation between the $(1,2)$- rainbow connection numbers of a graph and its spanning subgraph is shown in the following.
 \begin{proposition}[\cite{Li2016}] 
	Given a connected graph $G$. If $H$ is a connected spanning subgraph of $G$, then $rc_{1,l}(G) \leq rc_{1,l}(H)$. Particularly, $rc_{1,l}(G) \leq rc_{1,l}(T)$ in case $T$ is a spanning tree of $G$.
\end{proposition}

Combine Theorem \ref{theorem_tree}, Proposition \ref{p1} and the main result of Theorem \ref{thm1} and its corollaries, we have the followings.
 \begin{theorem}\label{thm2}
 	Let $G$ be a connected graph of order $n$ and $k\geq 2$ be an integer. If $\sigma_k(G)\geq n-1$, then $rc_{1,2}(G)\leq k+1.$ 
 \end{theorem}
 \section{Remarks}
 In the proof of Theorem \ref{thm1}, we may see that $\alpha(G)\geq k+1.$ This implies that if $\alpha(G)\leq k$ then $G$ has a spanning tree $T$ with $\Delta_2(T)\leq k+2.$ This also arises a natural question that how we can give a condition for $\sigma_{k+1}(G)$ to show that $G$ has $\Delta_2(G)\leq k+2.$ Thus, we propose the following problem.
 \begin{problem}\label{pro1} Let $G$ be a connected graph with order $n.$ If $\sigma_{k+1}(G)\geq \dfrac{k+1}{k}(n-2)+1,$ then $\Delta_2(G)\leq k+2.$
 \end{problem}
 We also note that if Problem \ref{pro1} is correct then the condition $\sigma_{k+1}(G)\geq \dfrac{k+1}{k}(n-2)+1$ is sharp by using the class of graph in Proposition \ref{p1}.
 
 On the other hand, the result of Theorem \ref{thm2} follows directly from Theorem \ref{thm1} and properties of generalized rainbow connection and, in Proposition \ref{p1}, we show a class of graph such that the result of Theorem \ref{thm1} is optimal. However, Zhu et al. \cite{Zhu2021} proved that the $(1,2)$-rainbow connection number of 2-connect graph is less then or equal $5.$ So the $(1,2)$-rainbow connection number of such graphs are bounded by $5.$ Therefore, it is natural to ask whether the condition “$\sigma_k(G)\geq n-1$" in Theorem \ref{thm2} is indeed optimal. If this is the case, what can be said about the corresponding classes of $G$?  In particular, we propose a question. 
\begin{problem} How can we tell about the graphs $G$ satisfying $\sigma_k(G)= n-2$ and $rc_{1,2}(G)\geq k+2$?
\end{problem}
  
Finally, let $n$, $k, k_0$, and $l$ represent arbitrary positive integers. We use $f(k)$ to denote a function of $k$, $g(n, k)$ for a function of both $n$ and $k$ and $h(k, m)$ for a function of both $k$ and $m$. We hope the subsequent problem can be addressed.
 
 \begin{problem}
 	\label{main_propblem}
 	Let $k\geq k_0, m\geq 3$. How we can find $f(k), g(n,k)$ and $h(k,m)$ such that: if $G$ is a connected graph of order $n$ and $\sigma_k(G)\geq g(n,k)$, then $rc_{1,m}G\leq h(k,m)$ for every  $n\geq f(k).$ 
 \end{problem}

\end{document}